\documentclass[12pt]{amsart}
\usepackage{amssymb,latexsym,amsmath,amsthm,amscd}
\usepackage{setspace}
\usepackage{verbatim}
\usepackage[active]{srcltx}
\usepackage[all]{xy} \xyoption{arc}
\usepackage[left=3cm,top=2cm,right=3cm,bottom = 2cm]{geometry}
\usepackage{graphicx}
\usepackage[usenames,dvipsnames]{color}
\usepackage{charter}

\theoremstyle{plain}

\theoremstyle{definition}

\theoremstyle{remark}

\newcommand{\CC}{\mathbf{C}}

\newcommand{\ZZ}{\mathbf{Z}}

\begin{document}
%Title
\title{A Brief History of the Positivity Conjecture in Tensor Category Theory}
\author{Geoffrey Mason}
\address{Department of Mathematics, UC Santa Cruz, CA 95064}
\email{gem@ucsc.edu}
\thanks{The author thanks the Simon Foundation, grant $\#427007$, for its support.}

\begin{abstract} We show the existence of a finite group $G$ having an irreducible character
$\chi$ with Frobenius-Schur indicator $\nu_2(\chi){=}{+}1$ such that $\chi^2$ has an irreducible constituent $\varphi$ with 
$\nu_2(\varphi){=}{-}1$.\ This provides counterexamples to the positivity conjecture in rational CFT and a conjecture of Zhenghan Wang about pivotal fusion categories.\\
MSC{:}$18D10, 20G05$
\end{abstract}
\maketitle

\section{Introduction}
There are many positivity conjectures in mathematics.\ The one that concerns us here
appeared at the turn of the century in  a paper of Peter Bantay (\cite{B}, eqn(2)).\ It states that
\begin{eqnarray}\label{posconj}
N_{pq}^r\nu_p\nu_q\nu_r{\geq}0.
\end{eqnarray}
The context is nominally that of rational CFT{:}\ $p, q, r$ are simple objects with Frobenius-Schur indicators
$\nu_p, \nu_q, \nu_r$ respectively, and $N_{pq}^r$ is the corresponding fusion rule.\ The conjecture is meaningful
in tensor categories equipped with appropriate additional structures, and we shall view it in that context.

\medskip
Soon after his paper appeared, Bantay brought the positivity conjecture to the author's attention.\ Framing it in a group-theoretic context, Bantay asked if one can find a finite group $G$ and a pair of irreducible characters $\alpha, \beta$ of $G$, each with Frobenius-SChur indicator ${+}1$, such that $\alpha\beta$ has an irreducible constituent  with indicator ${-}1$?\ I was able to supply such a $G$, thereby falsifying  (\ref{posconj}).\ The construction was  never published, but experts
became aware of its existence.

\medskip
In the meantime several variants of the positivity conjecture emerged.\ Already in 2001 Bantay had communicated a \emph{modified positivity conjecture}
to me, which he ascribed to `string theorists'.\ This conjecture proposed that (\ref{posconj}) holds
whenever $N_{pq}^r$ is \emph{odd}.\ This form of the  
conjecture is true and, as Richard Ng 
has pointed out,  may be deduced from the final Theorem in a 2002 paper by J\"{u}rgen Fuchs, Ingo Runkel and Christoph Schweigert  \cite{FRS}, who work in the context of a semisimple, braided, sovereign tensor category.

\medskip
More recently, Zhenghan Wang conjectured a version in his $2006$ book on topological
quantum computation  (\cite{W}, Conjecture 4.26)  which may be stated as follows (we are in a pivotal fusion category){:}
\begin{eqnarray*}
N_{p, p^{\vee}}^r>0 \Rightarrow \nu_r{=}1.
\end{eqnarray*}

\medskip
It transpires that the same counterexample to the positivity conjecture also serves as a counterexample to Wang's conjecture, and the purpose of this Note is  to belatedly make this example public.\ 
\medskip
We will show the following{:}
\begin{eqnarray*}
&&\mbox{There is a pair $(G, \chi)$ consisting of a finite group $G$ (of order $2^7$) }\\
&&\mbox{and an irreducible character $\chi$ of $G$ such that $\chi$ has indicator $\nu_2(\chi){=}{+}1$}\\ 
&&\mbox{and $\chi^2$ has an irreducible constituent $\varphi$ with  indicator $\nu_2(\varphi){=}{-}1$}.
\end{eqnarray*}

\medskip
We are indebted to Eric Rowell for bringing Wang's conjecture to our attention and Richard Ng 
for helpful correspondence.

\section{Construction} We use standard facts from representation theory, but include some details for the sake of completeness.

\medskip\noindent
Notation{:} if $G$ is a group and $a, b{\in}G$ then $[a, b]{:=}a^{-1}b^{-1}ab$.\\
$\hat{G}{:=}Hom(G, \CC^{\times})$ is the  group of \emph{characters} of $G$; $Z(G)$ is the \emph{center} of $G$.

\bigskip\noindent
1)\ Take $G$ to be any $2$-group with a normal subgroup
$H{\unlhd}G$ satisfying  $G/H{\cong}Q_8$, $H{\cong}\ZZ_2^4$ and $C_G(H){=}H$.\ For example,
we may take $G{=}H{\rtimes}Q$ where $Q{\cong}Q_8$ acts faithfully (by conjugation) on $H$.\ (The existence of
such a $G$ is shown below.)

\medskip\noindent
2)\ Let $\lambda{\in}\hat{H}$ be such that for each coset $xH{\not=}H$ in $G$ (alternatively, in the semi-direct product case, for $1{\not=}x{\in} Q$)
$\ker \lambda$ does \emph{not} contain $[H, x]$.\ (The existence of $\lambda$ is shown below. $[H, x]$ is defined below.)

\medskip\noindent
Claim 1.\ $\chi{:=}Ind_H^G\lambda$ is \emph{irreducible}.\\
Proof.\  One knows (\cite{CR}, Corollary 10.20) that
$\chi$ is irreducible, if, and only if, $\lambda{\not=} {^x}\lambda$ for all $x{\in}G{\setminus}{H}$.\ Here,
${^x}\lambda$ is defined by
\begin{eqnarray*}
{^x}\lambda(h){:=}\lambda(x^{-1}hx)\ \ (h{\in}H).
\end{eqnarray*}
Now we have
\begin{eqnarray*}
&&\lambda {=}{^x}\lambda \Leftrightarrow \lambda(h){=} {^x}\lambda(h)\ \ (h{\in}H)\\
&&\Leftrightarrow \lambda(h){=}\lambda(x^{-1}hx)\Leftrightarrow \lambda(h^{-1}x^{-1}hx){=}1\\
&&\Leftrightarrow [H, x]{\subseteq}ker\lambda,
\end{eqnarray*}
where
$[H, x]{:=}\langle [h, x]{\mid}h{\in}H\rangle$.

\medskip
By 2) this never happens.\ So indeed $\chi$ is irreducible. $\hfill\Box$

\medskip\noindent
Claim 2.\ $\chi^2$ has an irreducible constituent $\varphi$ with $\nu_2(\varphi){=}{-}1$.\\
Proof.\ Mackey has made a general study of the decomposition of tensor squares of induced characters such as
$\chi$ (\cite{CR}, Theorem 12.17).\ In particular,  it is known (loc.\ cit.)\ that there is a constituent equal to
$(\lambda^2)^G$.\  Now $\lambda$ is a character of $H$, which is a group of exponent $2$.\ Therefore, $\lambda^2{=}1_H$.\ It follows that $(\lambda^2)^G{=}(1_H)^G{=}\rho_{G/H}{=}\rho_{Q}$ is the \emph{regular representation} of
$G/H$.\ This has every irreducible character of $Q_8$ as a constituent, and in particular it has a constituent  with FS indicator equal to ${-}1$, namely the 2-dimensional irreducible $\varphi$ for $Q_8$ (lifted to $G$). $\hfill\Box$

\medskip\noindent
Claim 3.\ $G$ exists.\\
We need to show that $Q_8{\subseteq}GL_4(2)$.\ It is
enough to prove that $Q_8{\subseteq}A_8$(the alternating group of degree $8$), on account of the exceptional isomorphism $GL_4(2){\cong}A_8$.\ Consider the left regular representation of
$Q_8$. It provides an embedding $Q_8{\subseteq}S_8$.\ In fact, each element of order $4$ in $Q_8$
is represented as a product of a pair of cycles of length $4$, and such a permutation is \emph{even}, hence lies in $A_8$.\ 
Since $Q_8$ is generated by its elements of order $4$ we obtain $Q_8{\subseteq} A_8$. $\hfill\Box$

\medskip\noindent
Claim 4.\ 
 Let $1{\not=}z{\in}Z(Q)$, and set $H_0{:=}[H, z]$.\ Then $|H_0|{=}2$.\\
  This is the same as
 the assertion $|C_H(z)|{=}8$.\ To see this, let $P{\subseteq}GL_4(2)$ be the stabilizer of some  hyperplane of $H$, call it $H_1$.\ Then $P{\cong}D{\rtimes}L$ with $D{\cong}\ZZ_2^3, L{\cong}GL_3(2)$, moreover $D$ acts trivially on $H_1$.\ Notice that $P$ contains a Sylow $2$-subgroup of $GL_4(2)$, so that by Sylow's theorem we may assume that $Q{\subseteq}P$.\
 If $z{\in}D$ then our assertion follows, whereas if $z{\notin}D$ then $Q{\cap}D{=}1$, so that $Q$ is isomorphic to a Sylow $2$-subgroup of $P/D{\cong}L$.\ But this is impossible because a Sylow $2$-subgroup of $L_3(2)$ is dihedral of order $8$.
 $\hfill\Box$
 
 \medskip\noindent
 Claim 5.\ $\lambda$ exists.\\
  We assert that
 \begin{eqnarray}\label{inter}
H_0{=}\cap_{1\not=x\in Q}[H, x].
\end{eqnarray}
It suffices to show that $H_0{\subseteq}[H, x]$ whenever $x{\in}Q$ has order $4$.
To see this, suppose that $x$ has order $4$.\ Then $x^2{=}z$.\ We have
\begin{eqnarray*}
[h, x^2]{=}h^{-1}x^{-2}hx^2{=}(h^{-1}x^{-1}hx)(x^{-1}h^{-1}x^{-1}hx^2){=}[h, x][h, x]^x{=}[[h, x], x]{\in}[H, x],
\end{eqnarray*}
as needed. 

\medskip
To get the existence of $\lambda$, choose any hyperplane $H_2$ of $H$ that does \emph{not} contain $H_0$.\ $H_2$ exists thanks to Claim 4.\ Choose $\lambda{\in}\hat{H}$ so that
$H_2{:=}ker\lambda$.\ Since $H_0{\subseteq}[H, x]$
for all $x{\in}G{\setminus}{H}$ by (\ref{inter}) then $[H, x]{\not\subseteq} ker\lambda$. $\hfill\Box$

\medskip\noindent
Claim 6.\ $\nu_2(\chi){=}{+}1$.\\
We have $\nu_2(\chi){=}|G|^{-1}\sum_{g\in G}\chi(g^2)$.\ Because $\chi$ is induced from the normal subgroup $H$ then $\chi$ \emph{vanishes} on $G{\setminus}{H}$.\ Moreover, $g^2{\in}H$ if, and only if, $g{\in}H\langle z\rangle$.\ If
$g{\in}H$ then $\chi(g^2){=}\chi(1){=}8$.\ If $g{\in}Hz$, say $g{=}hz$, then $g^2{=}[h,z]{\in}H_0$.\ Moreover, there are just 8 choices of $h$ with $[h, z]{=}1$, and for these we again have $\chi(g^2){=}8$.\ As for the $8$ choices satisfying $g^2{=}[h, z]{\not=}1$, by construction we have $\lambda(g^2){=}{-}1$, moreover $g^2{\in}Z(G)$.\ Therefore\\ 
\begin{eqnarray*}
\chi(g^2){=}|H|^{-1}\sum_{t\in G} \lambda(t^{-1}g^2t){=}|H|^{-1}|G|\lambda(g^2){=}{-}8.
\end{eqnarray*}
The Claim follows from these calculations. $\hfill\Box$

\end{document}